\documentclass[12pt]{article}
\usepackage{amssymb, amsfonts, amsmath, amsthm, url, graphicx}

\def\1{^{-1}}
\def\barr{\begin{array}}
\def\earr{\end{array}}

\def\<{\left<}
\def\>{\right>}

\def\3{\subset }
\def\4{\subseteq }
\def\ov{\overline}

\def\0{\leqno}

\def\barr{\begin{array}}
\def\earr{\end{array}}

\def\nn{{\rm I\!N}}
\def\Z{{\rlap{$\kern2pt{\rm Z}$}{\rm Z}\,}}

\title{\bf On the converse\\ of Fuzzy Lagrange's Theorem}
\author{Marius T\u arn\u auceanu}
\date{February 17, 2015}

\begin{document}

\maketitle

\begin{abstract}
In fuzzy group theory many versions of the well-known Lagrange's
theorem have been studied. The aim of this article is to investigate
the converse of one of those results. This leads to an interesting
cha\-rac\-te\-ri\-za\-tion of finite cyclic groups.
\end{abstract}
\bigskip


\bigskip

\section{Introduction}

The notion of fuzzy set was introduced by Zadeh in 1965. The importance of the introduced notion of fuzzy set was realized by the research worker in all the branches of science and technology and has successfully been exploited. Recently fuzzy set theory has been applied in analysis and topology by Tripathy and Baruah \cite{10}, Tripathy and Borgohain \cite{11}, Tripathy and Das \cite{12}, Tripathy and Sarma \cite{13}, Tripathy, Baruah, Et and Gungor \cite{14}, Tripathy and Ray \cite{15} and many others.
\bigskip

One of the most important results of finite group theory is the
Lagrange's theorem. The fuzzification of this theorem has been
studied by several authors and a lot of results one could call a
Fuzzy Lagrange's Theorem have been obtained (see Section 2.3 of
\cite{7}). In the following we will focus on one of them, Theorem
2.3.17 of \cite{7}, that will be called the Fuzzy Lagrange's
Theorem.

\bigskip\noindent{\bf Fuzzy Lagrange's Theorem.} {\it Let $G$ be a finite group of order $n$.
Then $O(\mu)$ is a divisor of $n$, for every fuzzy subgroup $\mu$
of $G$.}
\bigskip

On the other hand, it is well-known that the converse of the
classical Lagrange's theorem is not true for all finite groups.
More precisely, the finite groups satisfying this, usually called
CLT-groups, determine a class between supersolvable groups and
solvable groups (see, for example \cite{1} and \cite{2}). In the
current paper we will investigate the analogue problem for fuzzy
subgroups. First of all, we formulate the converse of the above
theorem.

\bigskip\noindent{\bf Converse of Fuzzy Lagrange's Theorem.} {\it Let $G$ be a finite group of order $n$.
Then, for every divisor $d$ of $n$, there is a fuzzy subgroup
$\mu_d$ of $G$ such that $O(\mu_d)=d$.}
\bigskip

By taking $G$ of small order, we observe that the Converse of
Fuzzy Lagrange's Theorem is true for the cyclic groups
$\mathbb{Z}_2$, $\mathbb{Z}_3$, ... and so on, but it fails for
the Klein's group $\mathbb{Z}_2\times\mathbb{Z}_2$ (that has no
fuzzy subgroup of order 4) or for the symmetric group $S_3$ (that
has no fuzzy subgroup of order 3). In this way, the following
question is natural: \textit{which are the finite groups $G$
satisfying the Converse of Fuzzy Lagrange's Theorem}?
\bigskip

Our following main result completely answers this question, by proving that:

\bigskip\noindent{\bf Theorem.} {\it A finite group $G$ satisfies
the Converse of Fuzzy Lagrange's Theorem if and only if it is
cyclic.}
\bigskip

Remark that one obtains a new characterization of the finite
cyclic groups by using "fuzzy group ingredients" (notice also that
other such characterizations can be found in our previous papers
\cite{8} and \cite{9}).

\section{Preliminaries}

Let $(G,\cdot,e)$ be a group ($e$ denotes the identity of $G)$ and
$\mu:G\to[0,1]$ be a fuzzy subset of $G$.  We say that $\mu$ is a
\textit{fuzzy subgroup} of $G$ if it satisfies the following two
conditions:
\begin{enumerate}
\item[a)] $\mu(xy)\ge\min\{\mu(x),\mu(y)\},$ for all $x,y\in G$;
\item[b)] $\mu(x^{-1})\ge\mu(x)$, for any $x\in G$.
\end{enumerate}
In this situation we have $\mu(x^{-1})=\mu(x)$, for any $x\in G$,
and $\mu(e){=}\sup\mu(G)$. If $\mu$ satisfies the supplementary
condition
$$\mu(xy)=\mu(yx),\mbox{ for all }x,y\in G,$$then it is called a {\it fuzzy normal subgroup} of $G$. As in the case of subgroups,
the sets $FL(G)$ and $FN(G)$ consisting of all fuzzy subgroups and
of all fuzzy normal subgroups of $G$ are lattices with respect to
fuzzy set inclusion, called the {\it fuzzy subgroup lattice} and
the {\it fuzzy normal subgroup lattice} of $G$, respectively.
\bigskip

For each $\alpha\in [0,1]$, we define the level subset
$$_\mu G_\alpha=\{x\in G\mid\mu(x)\ge\alpha\}.$$These subsets
allow us to characterize the fuzzy (normal) subgroups of $G$, in
the following manner: $\mu$ is a fuzzy (normal) subgroup of $G$ if
and only if its level subsets are (normal) subgroups in $G$.
\bigskip

The concept of fuzzy order of an element $x\in G$ relative to a
fuzzy subgroup $\mu\in FL(G)$ has been defined in \cite{6}, as
follows.
\bigskip

If there exists $n\in \mathbb{N^*}$ such that
$\mu(x^n)=\mu(e)$, then $x$ is said to be of \textit{finite fuzzy
order} with respect to $\mu$ and the least such positive integer
$n$ is called the \textit{fuzzy order} of $x$ with respect to
$\mu$ and written as $FO_{\mu}(x)$. If no such $n$ exists, $x$ is
said to be of \textit{infinite fuzzy order} with respect to $\mu$.
Clearly, in a finite group $G$ all elements have finite fuzzy
orders relative to any fuzzy subgroup of $G$. Under the above
hypotheses, we also have
$$FO_{\mu}(x)=o_H(x),$$where $H=\{a\in G \mid \mu(a)=\mu(e)\}\leq G$ and $o_H(x)$ denotes
the order of $x$ relative to $H$ (i.e. the smallest positive
integer $n$ such that $x^n\in H$, if there exists such a positive
integer). In particular, if $H$ is the trivial subgroup $\{e\}$ of
$G$, then
$$FO_{\mu}(x)=o(x),$$the (classical) order of $x$ in $G$.
\bigskip

Let $\mu\in FL(G)$. If there exists $n\in \mathbb{N^*}$ such that
$\mu(x^n)=\mu(e)$, for all $x\in G$, then the smallest such
positive integer is called the \textit{fuzzy order} of $\mu$ and
written as $O(\mu)$. If no such positive integer exists, then
$\mu$ is said to be of \textit{infinite fuzzy order}. It follows
immediately that if $G$ is a finite group, then $$O(\mu)={\rm
lcm}\{FO_{\mu}(x)\mid x\in G\},$$or equivalently
$$O(\mu)={\rm exp}_H(G),$$the exponent of $G$ relative to $H=\{a\in G \mid
\mu(a)=\mu(e)\}\leq G$ (i.e. the least common multiple of the
orders of all elements of $G$ relative to $H$).

\section{Proof of the main theorem}

First of all we prove two preliminary results those will be used in establishing the main result.

\bigskip\noindent{\bf Lemma 1.} {\it Let $G$ be a finite group which satisfies
the Converse of Fuzzy Lagrange's Theorem. Then all Sylow subgroups
of $G$ are cyclic.}

\bigskip\noindent{\bf Proof.} Let $n$ be the order of $G$. By our
hypothesis, there is a fuzzy subgroup $\mu_n$ of $G$ satisfying
$O(\mu_n)=n$. This means there is a subgroup $H$ of $G$
satisfying ${\rm exp}_H(G)=n$. Since ${\rm exp}_H(G)\mid{\rm
exp}(G)$, one obtains $${\rm exp}(G)=n.\0(1)$$

Let $n=p_1^{\alpha_1}p_2^{\alpha_2}...p_k^{\alpha_k}$ be the
decomposition of $n$ as a product of prime factors and
$i\in\{1,2,...,k\}$. By (1), we infer that there exists $a\in G$
whose order is divisible by $p_i^{\alpha_i}$, say
$o(a)=p_i^{\alpha_i}q$ for some $q\in \mathbb{N}^*$. Then
$o(a^q)=p_i^{\alpha_i}$, that is $G$ contains an element of order
$p_i^{\alpha_i}$. In other words, $G$ possesses a cyclic Sylow
$p_i$-subgroup, and therefore all Sylow $p_i$-subgroups of $G$ are
cyclic. This completes the proof. \hfill\rule{1,5mm}{1,5mm}

\bigskip\noindent{\bf Remark.} In group theory, the finite groups all of whose
Sylow subgroups are cyclic are usually called ZM-$groups$ (see
\cite{4} and \cite{5}). Such a group is of type
$${\rm ZM}(m,n,r)=\langle a, b \mid a^m = b^n = 1,
\hspace{1mm}b^{-1} a b = a^r\rangle,$$where the triple $(m,n,r)$
satisfies the conditions $${\rm gcd}(m,n) = {\rm gcd}(m, r-1) = 1
\quad \text{and} \quad r^n \equiv 1 \hspace{1mm}({\rm
mod}\hspace{1mm}m).$$The subgroups of ${\rm ZM}(m,n,r)$ have been
completely described in \cite{3}. Let
$$L=\left\{(m_1,n_1,s)\in\nn^3 \hspace{1mm}\mid\hspace{1mm} m_1\mid
m,\hspace{1mm} n_1\mid n,\hspace{1mm} s<m_1,\hspace{1mm} m_1\mid
s\frac{r^n-1}{r^{n_1}-1}\right\}.$$Then there is a bijection
between $L$ and the lattice of subgroups of ${\rm ZM}(m,n,r)$,
namely the function that maps a triple $(m_1,n_1,s)\in L$ into the
subgroup $H_{(m_1,n_1,s)}$ defined by
$$H_{(m_1,n_1,s)}=\bigcup_{k=1}^{\frac{n}{n_1}}\alpha(n_1,
s)^k\langle a^{m_1}\rangle=\langle a^{m_1},\alpha(n_1,
s)\rangle,$$where $\alpha(x, y)=b^xa^y$, for all $0\leq x<n$ and
$0\leq y<m$.

\bigskip\noindent{\bf Lemma 2.} {\it Let $G={\rm ZM}(m,n,r)$ be a {\rm ZM}-group
satisfying the Converse of Fuzzy Lagrange's Theorem. Then $G$ is
cyclic.}

\bigskip\noindent{\bf Proof.} By applying our hypothesis for
$d=n$, we infer that there is a subgroup $H=H_{(m_1,n_1,s)}$ of
$G$ such that $${\rm exp}_H(G)=m.\0(2)$$On the other hand, we have
$${\rm exp}_H(G)={\rm lcm}\{o_H(a), o_H(b)\}=o_H(a)o_H(b),$$where
$a$ and $b$ are the generators of $G$. Since $o_H(a)\mid m$,
$o_H(b)\mid n$ and ${\rm gcd}(m,n)=1$, it follows that we must
have $o_H(a)=m$ and $o_H(b)=1$. This leads to $m_1=m$, $n_1=1$ and
$s=0$, that is $$H=H_{(m,1,0)}=\langle b \rangle.$$

Let $x=aba^{-1}\in G$. It is easy to see that
$$o_H(x)\mid n,\0(3)$$because $o(x)=o(b)=n$. By (2), one obtains that $$o_H(x)\mid
m.\0(4)$$Obviously, the relations (3) and (4) imply $o_H(x)=1$, in
view of the condition ${\rm gcd}(m,n)=1$. In this way, we have
$x\in H$, say $x=b^k$ for some integer $k$. It results that
$a^{r-1}=b^{k-1}$ and therefore $r=k=1$. This shows that
$b^{-1}ab=a$, i.e. $ab=ba$. Hence $G$ is cyclic, as desired.
\hfill\rule{1,5mm}{1,5mm}
\bigskip

We are now able to prove our main result.

\bigskip\noindent{\bf Proof of the main theorem.} If a finite group $G$ satisfies
the Converse of Fuzzy Lagrange's Theorem, then it is cyclic by
Lemma 1 and Lemma 2.

Conversely, let $G=\langle a \rangle$ be a cyclic group of order
$n$ and $d$ be an arbitrary divisor of $n$. Take $H_d=\langle a^d
\rangle$. Then we can easily see that $${\rm exp}_{H_d}(G)={\rm
exp}(G/H_d)=|G/H_d|=d.$$Define the fuzzy subset $\mu_d:
G\longrightarrow [0,1]$ by
$$\mu_d(x)=\left\{\barr{ll}1,&x\in H_d\\ 0,&x\in G\setminus H_d\\\earr\right., \hspace{1mm}\forall \hspace{1mm}x\in
G.$$Then the corresponding level subsets $$_{\mu_{d}} G_1=H_d \mbox{
and } _{\mu_d} G_0=G$$are subgroups of $G$, that is $\mu_d\in
FL(G)$. Moreover, we have $O(\mu_d)=d$. In other words, $G$
satisfies the Converse of Fuzzy Lagrange's Theorem.
\hfill\rule{1,5mm}{1,5mm}

\bigskip\noindent{\bf Remark.} It is well-known that all subgroups of a cyclic group are
normal. Then the above fuzzy subgroup $\mu_d$ is in fact a
fuzzy normal subgroup of $G$. This shows that the main theorem can
be reformulated in the following way: \textit{A finite group $G$
of order $n$ is cyclic if and only if\,, for every divisor $d$ of
$n$, there is a fuzzy normal subgroup $\mu_d$ of $G$ such that
$O(\mu_d)=d$}.
\bigskip

We note that our result is useful to describe
the fuzzy subgroups of finite cyclic groups.

\bigskip\noindent{\bf Example.} Let $\mathbb{Z}_{12}=\{\ov{0},\ov{1},...,\ov{11}\}$
be the additive group of integers modulo 12. Then, by the Fuzzy Lagrange's Theorem, we know that $O(\mu)\hspace{-0,5mm}\mid\hspace{-0,5mm}12$, $\forall\,\mu\in FL(\mathbb{Z}_{12})$. Since $\mathbb{Z}_{12}$ is a cyclic group, it also satisfies the Converse of Fuzzy Lagrange's Theorem: for every divisor $d$ of 12, there is $\mu_d\in FL(\mathbb{Z}_{12})$ defined as above such that $O(\mu_d)=d$. For example, a fuzzy subgroup of fuzzy order 3 of $\mathbb{Z}_{12}$ is\newpage $$\mu_3: \mathbb{Z}_{12}\longrightarrow [0,1],\,\,\, \mu_3(\ov{x})=\left\{\barr{ll}1,&\ov{x}\in \{\ov{0},\ov{3},\ov{6},\ov{9}\}\\ 0,&\ov{x}\in \{\ov{1},\ov{2},\ov{4},\ov{5},\ov{7},\ov{8},\ov{10},\ov{11}\}\\\earr\right..$$
\smallskip

Finally, we recall that the subgroup lattice of a finite cyclic group $G$ of order $n$ is isomorphic to the lattice $L_n$ of all divisors of $n$. A fuzzy version of this result can be also obtained by using our result, namely: $$\textit{The fuzzy subgroup lattice of $G$ contains a sublattice isomorphic to $L_n$}.$$Indeed, let $L=\{\mu_d\mid d\in L_n\}\subseteq FL(G)$. Then it is easy to see that $L$ is a sublattice of $FL(G)$. Moreover, the map $O:L\longrightarrow L_n$, $O(\mu_d)=d, \forall\, \mu_d\in L$, is a lattice isomorphism.

\section{Conclusions and further research}

The finite cyclic groups constitute one of the most famous classes
of finite groups. A large number of characterizations of these
groups is known (see, for example \cite{4} and \cite{5}). Our
main theorem gives another such cha\-rac\-te\-ri\-za\-tion, which
is based on "fuzzy group ingredients". It illustrates the
powerful connection between fuzzy group theory and group theory,
that is still developed in many works.
\bigskip

Remark also that the class of finite groups satisfying the
Converse of Fuzzy Lagrange's Theorem (namely the cyclic groups) is
different from the class of finite groups satisfying the Converse
of Lagrange's Theorem (namely the CLT-groups). So, the
characterization of finite cyclic groups given by this note is a
specific property of fuzzy group theory.
\bigskip

We end our paper by indicating an open problem concerning the fuzzy orders
of the fuzzy subgroups of a finite group.

\bigskip\noindent{\bf Open problem.} Let $G$ be a finite group of order $n$,
$L_n$ be the lattice of all divisors of $n$ and $O:
FL(G)\longrightarrow L_n$ be the map defined by $\mu\mapsto
O(\mu)$, $\forall \hspace{1mm}\mu\in FL(G)$. We observe that our
main theorem describes in fact the finite groups such that $O$ is
onto. What can be said in general about $Im(O)$? Does form it a
sublattice of $L_n$? Study other basic properties of this map, such as
injectivity, monotony, ... and so on.
\bigskip

\bigskip\noindent{\bf Acknowledgements.} The author is grateful to the reviewers for
their remarks which improve the previous version of the paper.

\vspace*{5ex}\small

\hfill
\begin{minipage}[t]{5cm}
Marius T\u arn\u auceanu \\
Faculty of  Mathematics \\
``Al.I. Cuza'' University \\
Ia\c si, Romania \\
e-mail: {\tt tarnauc@uaic.ro}
\end{minipage}

\end{document}